\documentclass[a4paper,12pt]{article}
\usepackage{amsfonts,amssymb, pdfpages,xcolor}

\title{Weak asymptotic methods for 3-D self-gravitating pressureless fluids. Application to the creation and evolution of solar systems from the fully nonlinear Euler-Poisson equations.}
\author{M. Colombeau,\\ \textit{mcolombeau@ime.usp.br,}\\
 Instituto de Matem\'atica e Estat\'istica,\\Universidade  de S\~ao Paulo, Brazil.}
\date{}
\begin{document}
\maketitle

\begin {abstract}
We construct a family of classical continuous functions $S(x,y,z,t,\epsilon)$ which tend to satisfy asymptotically the system of selfgravitating pressureless fluids when $\epsilon\rightarrow 0$. This produces a weak asymptotic method in the sense of Danilov, Omel'yanov and Shelkovich.  The construction  is based on   a family of  two ODEs (one for the continuity equation, one for the Euler equation) in  classical Banach spaces of continuous functions.    This  construction  applies to 3-D self-gravitating pressureless fluids even in presence of point and string concentrations of matter. The method is constructive which permits to check numerically from standard methods for ODEs that these functions tend to the known or admitted solutions when the latter exist. As a direct application we present a simulation of formation and evolution of a planetary system from a rotating disk of dust: a theorem in this paper asserts that the observed results are a depiction of functions that satisfy the system with arbitrary precision. 

\end{abstract}
This research was supported by  FAPESP, processo 2012/15780-9.\\

\textbf{1.  Introduction}.
 In this paper we consider selfgravitating fluids without pressure before the more delicate study in presence of pressure in \cite{Colombeaupressure}.
 A weak asymptotic method  is a sequence of approximate solutions which obey equations in the sense of distributions up to a small discrepancy that  tends to zero. Defined  by Danilov, Omel'yanov and Shelkovich as a continuation of  Maslov-Whitham asymptotic analysis   \cite{Danilov1}, they have recently been successfully used  by numerous authors \cite {Albeverio, Danilov1,Mitrovic, Shelkovich2, Shelkovich3,Omel'yanov,Panov, ShelkovichRMS,Shelkovichmat,Shelkovich1}, in particular for the study of creation and superposition of delta shocks arising in solutions of various systems. In absence of a uniqueness result the pertinence of the approximate solutions thus obtained is proved from the observations  that they  give the known  solutions.\\

   In this paper we construct weak asymptotic methods with full mathematical proofs for    systems of pressureless fluid dynamics on the $n$ dimensional torus $\mathbb{T}^n=\mathbb{R}^n/(2\pi \mathbb{Z})^n, n=1,2,3$ from an original method  stemming from the preliminary   numerical investigations in \cite{ColombeauSiam, ColombeauNMPDE, Colombeauideal}. Modifications in proofs give a similar result on the whole space $\mathbb{R}^n$ with initial conditions of finite total mass.
This  original method is based on the study of a family of two particular ODEs in the Banach space $\mathcal{C}(\mathbb{T}^n)$  of bounded continuous functions on $\mathbb{T}^n$ and suitable subspaces.  We obtain a family of continuous functions (of class $\mathcal{C}^1$ in $t$) $(x,y,z,t,\epsilon)\longmapsto S(x,y,z,t,\epsilon)$ which tend to satisfy the equations (in the distributional sense for the space variables) when $\epsilon\rightarrow 0^+$. For small fixed $\epsilon>0$ the functions $S(.,\epsilon)$ play therefore a role of approximate solutions. Furthermore these approximate solutions can be easily computed numerically from standard numerical schemes for ODEs. One observes the  coincidence with previously known solutions in all tests.\\

We first provide a weak asymptotic method for the system of pressureless fluids

\begin {equation}\frac{\partial \rho}{\partial t}+\vec{\nabla}.(\rho \vec{u})=0,\end{equation}
\begin {equation}\frac{\partial}{\partial t}(\rho \vec{u})+\vec{\nabla}.(\rho \vec{u} \otimes \vec{u})=\vec{0},\end{equation}
 where $\rho, \vec{u}=(u,v,w)$ denote respectively the density and the velocity vector.\\

 The result  extends  at once to the systems including a  pressureless energy equation \cite{Albeverio,  Nilsson1, Nilsson2, Shelkovich1}, i.e. (1,2) complemented by the pressureless energy equation 
  
\begin {equation}\frac{\partial }{\partial t}(\rho e)+\vec{\nabla}.(\rho e \vec{u})=0,\end{equation}
where $e$ denotes the total energy per unit mass. Information about its physical relevance  can be found in  \cite{Albeverio}. \\

 We obtain a weak asymptotic method as a family $( \rho^\epsilon,\vec{u^\epsilon},e^\epsilon)_\epsilon$  such that, for equation (1),  $\forall \psi\in \mathcal{C}_c^\infty(\mathbb{R}^3) , \forall t\in [0,+\infty[$,
 \begin {equation}\int  (\frac{\partial \rho^\epsilon}{\partial t}\psi+(\rho^\epsilon \vec{u^\epsilon}).\vec{\nabla}\psi) dxdydz\rightarrow 0,\end{equation}
 when $\epsilon\rightarrow 0^+$, and  other equations (2, 3). \\

 This result is extended to the system of self-gravitating pressureless fluids \cite{Coles} p. 207,  \cite{Peacock} p. 460, \cite{Peter}  p. 231, \cite{Charru}  p. 49:

\begin{equation}\frac{\partial \rho}{\partial t}+\vec{\nabla}.(\rho \vec{u})=0,\end{equation} 
\begin{equation}\frac{\partial}{\partial t}(\rho \vec{u})+\vec{\nabla}.(\rho \vec{u} \otimes \vec{u})
+\rho \vec{\nabla \Phi}=\vec{0},\end{equation}
\begin {equation}\Delta \Phi=4\pi G \rho, \end{equation}
where $\rho, \vec{u}=(u,v,w), \Phi$ denote respectively the density, the velocity vector and the gravitation potential;  $G$ is the universal gravitation constant.   These equations  are the continuity equation (5), the Euler equation (6) and the Poisson equation (7). The Poisson equation (Newton's law of gravitation) is complemented by boundary conditions which are  null at infinity in the planetary system simulation given below.   
 For the system (5-7) of self-gravitating pressureless fluids in several space dimensions  one obtains  a weak asymptotic method in 2-D even in presence of point concentrations of matter and in 3-D even in presence of point and string concentrations of matter.\\
 
 System (5-7) has already been considered by mathematicians and numerical physicists.
In 1-D the system of self-gravitating pressureless fluids has  been considered by mathematicians in \cite{Rykov}  and  \cite{Nguyen}  from a theoretical viewpoint. These authors have obtained results of existence of solutions under various assumptions.
 From the numerical viewpoint, cosmologists \cite{Coles, Peacock,Peter} and astrophysicists \cite{Heggie} have developped $N$-body simulations representing a sample of the universe as a box with periodic boundary conditions containing a large number of point masses interacting through their mutual gravity \cite{Coles}  pp. 304-310,  \cite{Peacock}  pp. 482-494. There exists a number of numerical codes for a self-gravitating collection of bodies. They represent a cosmological fluid as a discrete set of a large number of particles and  calculate the gravitational forces between them \cite{Blelloch}.  In absence of exact solutions for their validation, and impossibility of physical experiments, faith in these methods comes only from the fact they mimick the real physical process and reproduce qualitatively the aspect of the universe as it is observed, \cite{Coles}  p. 308.\\

The weak asymptotic method presented for system (5-7) extends in a staightforward way to the more general Euler-Poisson system in the case of several species of matter (dark matter and baryonic matter  in expanding background), as needed for the simulation of large structure formation in cosmology and for the evolution of galaxies. This weak asymptotic method has been  extended in presence of pressure to the 3-D systems of isothermal and isentropic fluids, and to the 2-D shallow water equations in \cite{Colombeaupressure}.\\


    Classical methods for the numerical solution of the ODEs under consideration, such as the explicit order one Euler method \cite{ColombeauSiam,ColombeauNMPDE} or, better, the RK4 Runge-Kutta method,  have given  exactly the known results \cite{Bouchut, Chertock, LeVeque} even in presence of point accumulation of matter or void regions.  This could be expected since the numerical method in \cite{ColombeauSiam, ColombeauNMPDE, Colombeauideal} was at the origin of the  method presented in this paper.\\
    
     As an application, from system (5-7) in 2-D, we present a numerical simulation of formation  of a planetary system from a rotating disk of dust: the rotating disk of dust collapses in a "star", concentrating the largest part of the matter from the disk and a finite number of "planets" that rotate around the star endlessly  with  variable trajectories   depending on the gravitational interaction between all these objects.\\


\textbf{2. A 1-D asymptotic formulation of the system of pressureless fluids.} In this section we state a family of systems of two ODEs whose solutions will provide a weak asymptotic method for system (1, 2) in 1-D, and we explain its origin from a natural approximation of the space derivative.\\

 For a map $u:\mathbb{T}^n\times \mathbb{R}^+\longmapsto\mathbb{R}$ we use the standard auxiliary notation
\begin{equation}u^+(x,t)=max(0,u(x,t)), \ u^-(x,t)=max(0,-u(x,t)),\end{equation}
then
\begin{equation}u(x,t)=u^+(x,t)- u^-(x,t), \ |u(x,t)|=u^+(x,t)+u^-(x,t).\end{equation}

We approximate the 1-D pressureless fluid system by the following systems of autonomous ODEs depending on a parameter $\epsilon$ when $\epsilon\rightarrow 0$, for which we will prove that $\rho(x,t,\epsilon)>0$,  thus permitting division by $\rho$:

\begin{equation}\frac{d}{dt}\rho(x,t,\epsilon)=\frac{1}{\epsilon}[(\rho u^+)(x-\epsilon,t,\epsilon)-(\rho |u|)(x,t,\epsilon)+(\rho u^-)(x+\epsilon,t,\epsilon)],\end{equation}

\begin{equation}\frac{d}{dt}(\rho u)(x,t,\epsilon)=\frac{1}{\epsilon}[(\rho u u^+)(x-\epsilon,t,\epsilon)-(\rho u |u|)(x,t,\epsilon)+(\rho u  u^-)(x+\epsilon,t,\epsilon)],\end{equation}

\begin{equation} u(x,t,\epsilon)=\frac{(\rho u)(x,t,\epsilon)}{\rho(x,t,\epsilon)}.\end{equation}

\textit{Remark}. In the particular case $u(x,t,\epsilon)$ has a constant sign, formulas (10, 11) reduce at once to the classical discretization of the $x-$derivative to the left for positive velocity and to the right for negative velocity. Formulas  (10, 11) will be justified a posteriori  by the result that they provide a weak asymptotic method for system (1, 2) in 1-D. Formulas (10, 11) are issued  from the scheme in \cite{ColombeauSiam} and therefore they can be obtained intuitively as follows. Consider cells of length $\epsilon$ centered respectively at the points $x,x-\epsilon$ and $x+\epsilon$. Let $\omega$ denote $\rho$ or $\rho u$ which are densities of mass and momentum respectively. The quantity $\epsilon.\omega(x,t+dt,\epsilon)$ in the cell of center $x$ at time $t+dt$ is equal to the quantity $\epsilon.\omega(x,t,\epsilon)$ at time $t$  minus the quantity $(\omega. |u|dt)(x,t,\epsilon)$ escaped from this cell between $t$ and $t+dt$ (to the left if $u<0$, to the right if $u>0$) plus the  quantities received by this cell from the left and from the right, that is $(\omega.u^+dt) (x-\epsilon,t,\epsilon), (\omega.u^-dt) (x+\epsilon,t,\epsilon)$ respectively, under the assumption $|u| dt \leq \epsilon$, so that the balance of matter in a cell involves the two closest neighbor cells only. Then one lets $dt$ tend to 0 for fixed $\epsilon$ to obtain (10, 11) with  $\rho$ and $\rho u$ in place of $\omega$. The above reasoning done in 3-D is identical to the reasoning of physicists to obtain (1, 2).\\

We approximate the initial conditions $\rho_0\in L^1(\mathbb{T}), \rho_0\geq 0$ and $ u_0\in  L^\infty(\mathbb{T})$ in their respective norms by a family $\{\rho_{0,\epsilon}, u_{0,\epsilon}\}_\epsilon$ of continuous  functions on $\mathbb{T}$ with the natural property $\rho_{0,\epsilon}(x)\geq \epsilon   \  \forall x$. \\

\textbf{3. A priori inequalities.}  In this section we derive the a priori inequalities that will permit to prove existence of a global flow for the solutions of  equations (10-12) in positive time. We denote by $\mathcal{C}(\mathbb{T})$ the Banach space of all bounded continuous real valued functions on $\mathbb{T}$ with the sup. norm. For fixed $\epsilon>0$ we assume the existence of a solution
$$[0,\delta(\epsilon)[ \longmapsto  (\mathcal{C}(\mathbb{T}))^2$$
\begin{equation} \ \  \ \ \  \ \ \ \ \ \ \ \ \ \ \   t  \longmapsto[x\mapsto(\rho(x,t,\epsilon),(\rho u)(x,t,\epsilon))]\end{equation}
continuously differentiable  on $[0,\delta(\epsilon)[$ (having a right hand-side derivative at $t=0$) with the following properties 
on the solution 
\begin{equation} \exists m>0  \ / \  \rho(x,t,\epsilon) \geq m  \  \ \forall x\in \mathbb{T}  \ \forall t\in [0,\delta(\epsilon)[,\end{equation}
\begin{equation} \exists M>0  \ /  \ \|(x\longmapsto u(x,t,\epsilon)) \|_ {L^\infty(\mathbb{T})} \ and \  \|(x\longmapsto \rho(x,t,\epsilon)) \|_ {L^\infty(\mathbb{T})}\leq M \ \ \forall t\in [0,\delta(\epsilon)[,\end{equation}
where  $u$ is defined by (12). We denote the initial conditions $ (x\longmapsto\rho(x,0,\epsilon))$ and $ (x\longmapsto u(x,0,\epsilon)) $ by $\rho_{0,\epsilon}$ and $u_{0,\epsilon}$ respectively. For fixed $\epsilon$ we obtain the following a priori estimates on the solution, depending on the initial condition and on time; the point is that $m$ and $M$  disappear in these inequalities.\\

\textbf{Proposition 1}. $\forall t\in [0,\delta(\epsilon)[, \forall x \in \mathbb{T}$, \textit{one has:}
\begin{equation}  \ |u(x,t,\epsilon)|\leq \|u_{0,\epsilon}\|_\infty, \ \ \ \ \ \ \ \ \ \ \ \ \ \ \ \ \ \ \ \ \ \ \ \ \ \ \ \ \ \ \ \ \ \ \ \ \ \ \ \ \ \ \ \ \ \ \ \ \ \  \ \ \ \ \ \ \ \ \ \ \ \end{equation}
\begin{equation}  \ \rho_{0,\epsilon}(x) exp(-\frac{\|u_{0,\epsilon}\|_\infty}{\epsilon}t)\leq \rho(x,t,\epsilon)\leq\|\rho_{0,\epsilon}\|_\infty exp(\frac{2\|u_{0,\epsilon}\|_\infty}{\epsilon}t),\end{equation}
\begin{equation}  \ \| \rho(.,t,\epsilon)\|_{L^1(\mathbb{T})} \leq \|\rho_{0,\epsilon}\|_{L^1(\mathbb{T})}. \ \ \ \ \ \ \ \ \ \ \ \ \ \ \ \ \ \ \ \ \ \ \ \ \ \ \ \ \ \ \ \ \ \ \ \ \ \ \ \ \ \ \ \ \ \ \ \ \  \end{equation}\\

\textit{Proof.} From (10) and the continuous differentiability assumption (13) on the solution, one has for fixed $\epsilon$: \\

$\rho(x,t+dt,\epsilon)=\rho(x,t,\epsilon)+\frac{dt}{\epsilon}[(\rho u^+)(x-\epsilon,t,\epsilon)-(\rho |u|)(x,t,\epsilon)+(\rho u^-)(x+\epsilon,t,\epsilon)]+dt.o(x,t,\epsilon,dt),$\\
\\
 where $\|o(.,t,\epsilon,dt)\|_\infty \rightarrow 0$ when $dt\rightarrow 0$ uniformly for $t\in[0,\delta']$ with  $\delta'<\delta(\epsilon)$  
(recall $\epsilon$ is fixed). This follows from the mean value theorem: if $f$ is a $\mathcal{C}^1$ function then $\|f(t+dt)-f(t)-f'(t).dt\| \leq sup_{0<\theta<1}\|(f'(t+\theta dt)-f'(t)).dt\|$ and from the uniform continuity of $\frac{d\rho(.,t,\epsilon)}{dt}$ in $t$, valued in the Banach space $\mathcal{C}(\mathbb{T})$, when $t$ ranges in a compact interval  $[0,\delta']$. Therefore
 \\

$\rho(x,t+dt,\epsilon)=
\frac{dt}{\epsilon}(\rho u^+)(x-\epsilon,t,\epsilon)+(1-\frac{dt}{\epsilon}|u|(x,t,\epsilon))\rho(x,t,\epsilon)+ $
\begin{equation}  \frac{dt}{\epsilon}(\rho u^-)(x+\epsilon,t,\epsilon)+dt.o(x,t,\epsilon,dt).\end{equation}
 For $dt>0$ small enough (depending on $\epsilon$),  so that the strictly positive  term $[1-\frac{dt}{\epsilon}|u|(x,t,\epsilon)]\rho(x,t,\epsilon)$ dominates the small term $dt.|o(x,t,\epsilon,dt)|$. Using (14, 15), we obtain, by inverting (19),\\

$  \frac{1}{\rho(x,t+dt,\epsilon)}= [\frac{dt}{\epsilon}(\rho u^+)(x-\epsilon,t,\epsilon)+$
\begin{equation}(1-\frac{dt}{\epsilon}|u|(x,t,\epsilon))\rho(x,t,\epsilon)+\frac{dt}{\epsilon}(\rho u^-)(x+\epsilon,t,\epsilon)]^{-1}+dt.o_1(x,t,\epsilon,dt)\end{equation} 
where $\|o_1(.,t,\epsilon)(dt)\|_\infty \rightarrow 0$ when $dt\rightarrow 0$ uniformly for $t\in[0,\delta']$. Therefore, using (19) stated with $\rho u$ in place of $\rho$, and (15),\\

$\frac{\rho u}{\rho}(x,t+dt,\epsilon)=\frac{\frac{dt}{\epsilon}(\rho u u^+)(x-\epsilon,t,\epsilon)+[1-\frac{dt}{\epsilon}|u|(x,t,\epsilon)](\rho u)(x,t,\epsilon)+\frac{dt}{\epsilon}(\rho u u^-)(x+\epsilon,t,\epsilon)}{\frac{dt}{\epsilon}(\rho u^+)(x-\epsilon,t,\epsilon)+[1-\frac{dt}{\epsilon}|u|(x,t,\epsilon)]\rho(x,t,\epsilon)+\frac{dt}{\epsilon}(\rho u^-)(x+\epsilon,t,\epsilon)}$+\begin{equation}dt.o_2(x,t,\epsilon,dt),\end{equation}
where $o_2$ has the same property as $o$ and $o_1$ above. Since for $dt>0$ small enough the above quotient is a barycentric combination of $u(x-\epsilon,t,\epsilon),  u(x,t,\epsilon)$ and $ u(x+\epsilon,t,\epsilon)$ which are in numerator inside $\rho u$, it follows that 
 \begin{equation} \|u(.,t+dt,\epsilon)\|_\infty\leq  \|u(.,t,\epsilon)\|_\infty+dt.\|o_2(.,t,\epsilon,dt)\|_\infty.\end{equation}

\textit{Lemma 1}. Let $f:\mathbb{R}^+\longmapsto\mathbb{R}^+$ be a function such that
\begin{equation} f(t+dt)\leq f(t)+dt.o(t,dt)\end{equation}
where $o(t,dt)\rightarrow 0$ uniformly in $t$ when $dt\rightarrow 0$. Then 
\begin{equation} \forall \tau>0 \ \  f(t+\tau)\leq f(t).\end {equation}

\textit{Proof of the lemma}. $f(t+\tau)=f(t)+\sum_{i=1}^n[f(t+i\frac{\tau}{n})-f(t+(i-1)\frac{\tau}{n})].$ 
\\
 Therefore (23) implies  $f(t+\tau)\leq f(t)+n \frac{\tau}{n}{o}(t,\frac{\tau}{n}) $ where ${o}(t,\frac{\tau}{n})\rightarrow 0$ uniformly in $t$ when $n\rightarrow \infty$ 
from the uniformness of a bound of $o(t,dt)$ in $t \in [0,\delta'[$. The result is obtained by letting $n\rightarrow \infty$.$\Box$\\

Application of the lemma to (22) on $[0,\delta']$ and letting $\delta'$ tend to $\delta(\epsilon)$ yields  $\|u(.,t+\tau,\epsilon)\|_\infty\leq  \|u(.,t,\epsilon)\|_\infty  \  \ \forall t \ and \ t+\tau \in [0,\delta(\epsilon)[$; in particular 
 $\|u(.,t,\epsilon)\|_\infty\leq  \|u_{0,\epsilon}\|_\infty$ which proves assertion (16).\\

Now let us prove assertion (17). From (10) $$\frac{d}{dt}\rho(x,t,\epsilon)\geq -\frac{1}{\epsilon}(\rho |u|)(x,t,\epsilon)$$ since $\rho,u^+$ and $u^-$ are positive. Therefore, from (16)
\begin{equation}\frac{d}{dt}\rho(x,t,\epsilon)\geq -\frac{\|u_{0,\epsilon}\|_\infty}{\epsilon}\rho(x,t,\epsilon).\end{equation}
Let $v(x,t,\epsilon)=\rho(x,0,\epsilon)exp(-\frac{\|u_{0,\epsilon}\|_\infty}{\epsilon}t)$. Then, using assumption (14) to divide by $\rho$,\\
\begin{equation}\frac{\rho_t}{\rho}(x,t,\epsilon)\geq\frac{v_t}{v}(x,t,\epsilon)=-\frac{\|u_{0,\epsilon}\|_\infty}{\epsilon}.\end{equation}
 By integration, since $\rho$ and $v$ have same initial condition and are positive, $log(\rho)\geq log (v)$, i.e. $\rho(x,t,\epsilon)\geq v(x,t,\epsilon)$, i.e.
\begin{equation}\rho(x,t,\epsilon)\geq \rho(x,0,\epsilon)exp(-\frac{\|u_{0,\epsilon}\|_\infty}{\epsilon} t),\end{equation}
which is the left hand-side inequality  (17).\\

 Now let us  prove the right hand-side inequality  (17).
From the positiveness of the two terms $\rho u^{\pm}$ in (10) and from (16),  one has

\begin{equation} \rho(x,t,\epsilon)\leq \rho_0(x,\epsilon)+\frac{2}{\epsilon}\int_0^t\|\rho(.,s,\epsilon)\|_\infty \|u_{0,\epsilon}\|_\infty ds.   \end{equation}
Since $\rho_0(x,\epsilon) \leq\| \rho_0(.,\epsilon)\|_\infty$ this implies\\

$\rho(x,t,\epsilon)\leq \| \rho_0(.,\epsilon)\|_\infty+\frac{2}{\epsilon}\int_0^t\|\rho(.,s,\epsilon)\|_\infty \|u_{0,\epsilon}\|_\infty ds$.\\
\\
 Since this holds for all $x$
\begin{equation} \|\rho(.,t,\epsilon)\|_\infty\leq \| \rho_0(.,\epsilon)\|_\infty+\frac{2}{\epsilon}\|u_{0,\epsilon}\|_\infty \int_0^t\|\rho(.,s,\epsilon)\|_\infty ds.   \end{equation}
Gronwall's inequality implies
\begin{equation} \|\rho(.,t,\epsilon)\|_\infty\leq \| \rho_0(.,\epsilon)\|_\infty exp(\frac{2}{\epsilon}\|u_{0,\epsilon}\|_\infty t).\end{equation}
\\
 Then we prove the $L^1$ bound (18). \\

 $\frac{d}{dt}\int_{-\pi}^{+\pi}\rho(x,t,\epsilon)dx=\frac{1}{\epsilon}\int_{-\pi}^{+\pi}[(\rho u^+(x-\epsilon,t,\epsilon)-\rho u^+(x,t,\epsilon))+(-\rho u^-(x,t,\epsilon)+\rho u^-(x+\epsilon,t,\epsilon))]dx=0$ \\
\\
from periodicity, which concludes the proof of Proposition 1. $\Box$\\
\\

\textbf{4. Global existence-uniqueness result for fixed $\epsilon$.} We  prove that the a priori inequalities permit to obtain existence of a global flow of class $\mathcal{C}^0$ for system (10-12) from the classical theory of ODEs in Banach spaces. In this section we use the notations: $X(x,t,\epsilon):=\rho(x,t,\epsilon), \ Y(x,t,\epsilon):=(\rho u)(x,t,\epsilon), \ X':=\frac{\partial X}{\partial t}, \ Y':= \frac{\partial Y}{\partial t}$. We assume properties (14, 15) on the initial condition only, i.e.  for $t=0$. Since $u=\frac{Y}{X}$ as long as $X\not=0$, which will be the case here, 
equations (10-12) become
\begin{equation} X'(x,t,\epsilon)=\frac{1}{\epsilon}[(Y^+)(x-\epsilon,t,\epsilon)-(|Y|)(x,t,\epsilon)
+(Y^-)(x+\epsilon,t,\epsilon)]\end{equation}
\begin{equation} Y'(x,t,\epsilon)=\frac{1}{\epsilon}[(Y(\frac{Y}{X})^+)(x-\epsilon,t,\epsilon)-(Y|\frac{Y}{X}|)(x,t,\epsilon)
+(Y(\frac{Y}{X})^-)(x+\epsilon,t,\epsilon)].\end{equation}
 If $0<\lambda<1$ let $\Omega_\lambda$ be the open set in $\mathcal{C}(\mathbb{T})$ defined by
\begin{equation} \Omega_\lambda:=\{(X,Y)\in(\mathcal{C}(\mathbb{T}))^2 / \ \forall x\in \mathbb{T}, \ \ \lambda<X(x)<\frac{1}{\lambda},\ and \ |Y(x)|<\frac{1}{\lambda}\}.\end{equation}
From the assumptions on the initial conditions, i.e. (14, 15) with $t=0$, it follows  that for given $\epsilon>0$ there is   some $\lambda>0$ such that
\begin{equation} (\rho_{0,\epsilon},\rho_{0,\epsilon}u_{0,\epsilon})\in \Omega_\lambda.\end{equation}
For convenience equations (31, 32) are stated as the following autonomous system, where  $F$ and $G$, given by (31, 32), map $\cup_{0<\lambda<1}\Omega_\lambda$ into $\mathcal{C}(\mathbb{T})$:
\begin{equation}X'(t)=F(X(t),Y(t)),\end{equation}
\begin{equation}Y'(t)=G(X(t),Y(t)).\end{equation}
The functions $F$ and $G$ in second members of (35, 36) are of class $\mathcal{C}_0$ from the absolute values involved in $(\frac{Y}{X})^\pm$. If $E_1$ and $E_2$ are  normed spaces, and if $\Omega$ is an open set in $E_1$ and if $k>0$ we denote by $Lip_k(\Omega,E_2)$ the set of all maps $f:\Omega\longmapsto E_2$ such that $\|f(x)-f(y)\|_{E_2}\leq k\|x-y\|_{E_1} \ \forall x,y\in \Omega$. \\

\textit{Lemma 2}. \textit{The maps $F$ and $G$ are in $Lip_{\frac{1}{\lambda^3}}(\Omega_\lambda \times \Omega_\lambda,\mathcal{C}(\mathbb{T})), \ if \ 0<\lambda<1$}.\\

\textit{Proof.} The map $(X,Y)\in \Omega_\lambda \times \Omega_\lambda\longmapsto\frac{Y}{X}\in \mathcal{C}(\mathbb{T}) $  satisfies $ |\frac{Y_1}{X_1}-\frac{Y_2}{X_2}|\leq\frac{1}{\lambda^3}(|X_2-X_1|+|Y_2-Y_1|)$. The maps $X\longmapsto|X|,X^+,X^-$ and the translation maps are in $Lip_1(\Omega_\lambda , \mathcal{C}(\mathbb{T}))$.$\Box$\\

\textbf{Proposition 2: global solution and flow.} \textit{The system of autonomous ODEs (35, 36) with initial condition (34) at $t=0$ admits a unique global solution on $[0,+\infty[$,     valued in} $\mathcal{C}(\mathbb{T})\times\mathcal{C}(\mathbb{T}).$\\

\textit{Proof.} Various versions of the proof following directly from Proposition 1 can be found in the litterature of ODEs. From the  theory of ODEs in Banach spaces, in  the Lipschitz case, for fixed $\epsilon>0$, the system of ODEs (35, 36) with initial condition (34) at $t=0$, admits a unique local solution on some interval $[0,\delta[$ 
valued in  some $\Omega_{\lambda}$: this local solution satisfies bounds (14, 15) with $m=\lambda, M=\frac{1}{\lambda}$ (33), if $t\in [0,\delta[$.\\

  Since the iteration proof of the Lipschitz case gives a uniform lower bound for the time existence when the initial conditions lie in $\Omega_\nu$,  it follows  that for any $\nu>0$, any $t_0\geq 0$  and any initial condition  in $\Omega_\nu$ at time $t_0$, there are $\delta'>0$ and $0<\nu'<\nu$, depending only on $\nu$, not on the initial conditions in $\Omega_\nu$ and not on $t_0$, such that there is a solution on $[t_0,t_0+\delta'[$ valued in $\Omega_{\nu'}$.\\
  
    By absurd, let $(X_{max},Y_{max})$ be a maximal solution on some interval $[0,T_{max}[$,  for $T_{max}<\infty$. Then, the a priori inequalities (16,  17) prove the existence of some $\nu>0$ such that $(X_{max}(t),Y_{max}(t))\in \Omega_{\nu} \ \forall t\in [0,T_{max}[$. Therefore there exists $\delta'>0 $ and $0<\nu'<\nu $ such that the solution $ (X_{max},Y_{max})$ can be extended to some interval  $[0,T_{max}+\delta'[,$ valued in  $\Omega_{\nu'}$. 
 Finally one obtains a unique global solution on $[0,+\infty)$ which further satisfies the bounds  in Proposition 1. $\Box$\\

\textit{Remark.} System (10-12) is a  family of differentiable dynamical systems in the Banach space $\mathcal{C}(\mathbb{T})^2$ depending continuously on the parameter $\epsilon>0$. Since the weak asymptotic method makes sense only when $\epsilon\rightarrow 0$ it appears as a germ at $\epsilon=0^+$ of these dynamical systems. The investigation of possible bifurcations for fixed $t$ when $\epsilon\rightarrow 0$ and of possible attractors when $t\rightarrow +\infty$ as in \cite{Temam} appears to be difficult to study mathematically. \\

\textbf{5. Weak asymptotic method.}  In this section we check that the solution of the ODEs provide a weak asymptotic method. By definition of a weak asymptotic method \cite{Danilov1} one has to prove that for a given  $\psi\in \mathcal{C}_c^\infty(\mathbb{R})$, 
\begin{equation} \int \frac{d}{dt}\rho(x,t,\epsilon)\psi(x)dx=\int \rho u(x,t,\epsilon)\psi'(x)dx +f(\epsilon)\end{equation}
and 
\begin{equation}  \int \frac{d}{dt}(\rho u)(x,t,\epsilon)\psi(x)dx=\int \rho u^2(x,t,\epsilon)\psi'(x)dx +g(\epsilon)\end{equation}
where $f(\epsilon)$  and $ g(\epsilon)\rightarrow 0$ when $\epsilon\rightarrow 0$. Since $|u|=u^++u^-$, equation (10) can be rewritten by changes in the integration variables\\ 

$\int \frac{d}{dt}\rho(x,t,\epsilon)\psi(x)dx=\frac{1}{\epsilon}[ \int(\rho u^+)(x,t,\epsilon)(\psi(x+\epsilon)-\psi(x))dx- \int(\rho u^-)(x,t,\epsilon)(\psi(x)-\psi(x-\epsilon))dx]$.\\
\\
Since $\rho u=\rho u^+-\rho u^-$, $\rho |u|=\rho u^++\rho u^-$ and $\frac{\psi(x+\epsilon)-\psi(x)}{\epsilon}=\psi'(x)+O(\epsilon)$ one obtains\\

 $\int \frac{d}{dt}\rho(x,t,\epsilon)\psi(x)dx= \int(\rho u)(x,t,\epsilon)\psi'(x)dx+O(\epsilon)\int_{K}\rho |u|(x,t,\epsilon)dx$\\
\\
 for another $O(\epsilon)$, where $K$ is a compact interval containing the support of $\psi$ in its interior.\\
 \\
From the $L^1-$stability of $\rho u$ from (16, 18) we obtain (37). The same proof holds for (38).\\

The advantage on the result in \cite{ColombeauSiam} obtained from a numerical scheme instead of ODEs lies in that one has now a usual derivative in time:  only the space  derivative is considered in the sense of distributions.  The important point is that the ODE method will provide proofs not needing  boundedness of velocity, which will be needed in the 2-D and  3-D selfgravitating case below and in presence of pressure \cite{Colombeaupressure}, thus constructing for the first time weak asymptotic methods with full proofs for the Euler-Poisson system.\\


We consider initial conditions $\rho_0 \in L^1(\mathbb{T}), \  \rho_0$  positive ($\rho_0$ can take null values) and  $  u_0 \in L^\infty(\mathbb{T})$. We approximate them in the sense of distributions by a family $\{\rho_{0,\epsilon} \in \mathcal{C}(\mathbb{T}), u_{0,\epsilon}\in \mathcal{C}(\mathbb{T})\}_\epsilon$ with the properties

\begin{equation} \rho_{0,\epsilon}(x)\geq \epsilon  \ \forall x, \ \  \|\rho_0-\rho_{0,\epsilon}\|_{L^1(\mathbb{T})} \rightarrow 0\end{equation}
 and
\begin{equation} \|u_{0,\epsilon}-u_0\|_{L^\infty(\mathbb{T})} \rightarrow 0\end{equation}
 when $\epsilon\rightarrow 0$.\\

From (39) note that $\epsilon$ here plays the role of $m$ in (14) (this causes no problem for (16-17) since $m$ disappears there), which permits, at the limit $\epsilon\rightarrow 0$, to consider void regions in the initial condition, although the  proof of Proposition 1 does not allow void regions. In the same spirit note that the presence of $\epsilon$ in the denominator in the right handside of  (17) permits at the limit $\epsilon\rightarrow 0$ the presence of point   concentrations of matter in the solution. Summarizing results from the previous sections  we have obtained:\\

\textbf{Theorem 1.}\textit{ Under  the above assumptions  on the initial conditions the system of ODEs (10, 11) provides a weak asymptotic method of order 1 for  the system  of 1-D pressureless  fluids, i.e.(1, 2) restricted to 1-D,  
which is global in space and positive time.}\\

The proof applies also to the pressureless energy equation  (3).\\
\\

\textbf{6. Extension of the ODE method to 2 and 3 space dimension.}  In this section we extend the statement of the family of ODEs to 2-D and 3-D by analogy with the 1-D case and we show that we obtain similar results. As exposed  in  section 2 the second members of formulas (10, 11) can be understood as a balance of the quantities $\omega=\rho$ and $\omega=\rho u $ respectively: $(\omega u^+)(x-\epsilon)dt$ that come from the left,  $(\omega u^-)(x+\epsilon)dt$ that come from the right   and    $(\omega |u|)(x)dt$  that escape from the center during time $dt$,  
 in analogy with the numerical scheme in \cite{ColombeauSiam}. This intuitive description and analogy with the 2-D and 3-D scheme in \cite{ColombeauSiam} will permit to extend easily the second members of (10, 11) to  2-D and 3-D.\\

\textit{Extension of the ODE to 2-D.} Let $S$ be a small square of side $\epsilon$  of center $(x,y)$. Let us denote by $u$ and $v$ the components of the velocity in the directions $(Ox)$ and $(Oy)$ respectively. The square $S$ has 4 neighbor squares having an edge in common with it and 4 neighbor squares having a vertex in common with it. We assume $\|u\|_\infty dt\leq \epsilon$ and $\|v\|_\infty dt\leq \epsilon$ so that between $t$ and $t+dt$ the square $S$ cannot receive matter coming from the squares that are not these 8 direct neighbors. This assumption can be easily satisfied since we will show that the maximum principle holds as in 1-D above, which makes the above bounds on velocity satisfied as soon as they are satisfied in the initial conditions. The 2-D balance evaluation between times $t$ and $t+dt$ extending the 1-D remark in section 2 is:

\ \\
$\epsilon^2\omega(x,y,t+dt,\epsilon)= \epsilon^2\omega(x,y,t,\epsilon)- (amount \ of \  \omega \ that \ escapes \  from \  S)+$\\ 
\\
$(amount \ of \  \omega \ that\ comes  from \  the\ 4 \ edge  \ neighbors) +(amount \ of  \omega \   that $\\ \begin{equation} comes \ from \  the  \ 4 \ vertex  \ neighbors).\end{equation}

 Let us evaluate the various terms in the second member above when the side of the square $S$ has length $\epsilon$ and the duration time is assumed to be $dt$ (these values give (10, 11) in the 1-D case). This can be done easily with the aid of  pictures such as figure A1 in \cite{ColombeauNMPDE} p. 97 and it gives the following results:\\

a) $\omega$ that escapes from the square $S$: we obtain $\omega(x,y,t,\epsilon)[|u(x,y,t,\epsilon)|dt\epsilon+|v(x,y,t,\epsilon)|dt \epsilon-|u(x,y,t,\epsilon)||v(x,y,t,\epsilon)|(dt)^2]$. To check this formula it suffices to draw a  picture of the translation $\mathcal{T}_{udt,vdt}S$ of $S$ by the vector $(u(x,y,t,\epsilon) dt, v(x,y,t,\epsilon) dt)$ and evaluate the area of $\mathcal{T}_{udt,vdt}S \cap \mathcal{C}S$ where $\mathcal{C}S$ is the complement of $S$.\\

b) $\omega$ that comes from the 4 edge neighbor squares: we obtain $(\omega. u^+dt.(\epsilon-|v|dt))(x-\epsilon,y,t,\epsilon) +(\omega. u^-dt.(\epsilon-|v|dt))(x+\epsilon,y,t,\epsilon)+(\omega. v^+dt.(\epsilon-|u|dt))(x,y-\epsilon,t,\epsilon)+(\omega. v^-dt.(\epsilon-|u|dt))(x,y+\epsilon,t,\epsilon).$ To check this formula one evaluates the area of the intersection with $S$  of the translated of each of these neighbor squares by the vector $(u(x',y',t,\epsilon) dt, v(x',y',t,\epsilon) dt)$  of their own velocity multiplied by time $dt$ when $ (x',y')=(x-\epsilon,y), (x+\epsilon,y), (x,y-\epsilon) $ and $(x,y+\epsilon)$ successively. \\

c) $\omega$ that comes from the 4 vertex neighbors: similarly we obtain $(\omega.u^+dt. v^+dt)(x-\epsilon,y-\epsilon,t,\epsilon)+
(\omega. u^+dt. v^-dt)(x-\epsilon,y+\epsilon,t,\epsilon)+(\omega. u^- dt.v^-dt)(x+\epsilon,y+\epsilon,t,\epsilon)+(\omega. u^- dt.v^+dt)(x+\epsilon,y-\epsilon,t,\epsilon)$. \\

The 1-D case  (10, 11) is obtained by letting $v=0$:\\
\\
$\bullet$ the estimate a) gives 
$(\omega |u|)(x,t,\epsilon)dt \epsilon$;\\$\bullet$ the estimate b) gives $(\omega u^+(x-\epsilon,t,\epsilon)+\omega u^-(x+\epsilon,t,\epsilon))dt \epsilon$;\\$\bullet$ the estimate c) is irrelevant in the 1-D case.\\

Using the 2-D balance evaluation  with the terms a), b) and  c) above gives the following statement for $\omega=\rho, \rho $ and , $\rho v$:\\

$ \frac{\omega(x,y,t+dt,\epsilon)-\omega(x,y,t,\epsilon)}{dt}=\frac{1}{\epsilon}[-(\omega(|u|+|v|))(x,y,t,\epsilon)+(\omega u^+)(x-\epsilon,y,t,\epsilon)+ $\\
\\
$ (\omega u^-)(x+\epsilon,y,t,\epsilon)+ (\omega v^+)(x,y-\epsilon,t,\epsilon)+(\omega v^-)(x,y+\epsilon,t,\epsilon)] +$ \\
\\
 $ \frac {dt}{\epsilon^2}[\omega |u||v|)(x,y,t,\epsilon)
-(\omega u^+|v|)(x-\epsilon,y,t,\epsilon)-  (\omega u^-|v|)(x+\epsilon,y,t,\epsilon)- (\omega v^+|u|)(x,$\\
\\
$y-\epsilon,t,\epsilon)-(\omega v^-|u|)(x,y+\epsilon,t,\epsilon)+ (\omega u^+ v^+)(x-\epsilon,y-\epsilon,t,\epsilon)+(\omega u^+ v^-)(x-$\begin{equation} \epsilon,y+\epsilon,t,\epsilon)+
 (\omega u^- v^-)(x+\epsilon,y+\epsilon,t,\epsilon)+ (\omega u^- v^+)(x+\epsilon,y-\epsilon,t,\epsilon)].   \end{equation}

 Formula (42) has been derived under the assumption $dt\|u\|_\infty\leq \epsilon$ and $ dt\|v\|_\infty\leq \epsilon$ needed for the physical interpretation that, in time duration $dt$, transport occurs  only between one cell and its 8  neighbor cells (having in common  an edge or a vertex). Formula (42) is directly inspired from \cite{ColombeauSiam} (22-24).
To obtain the differential equation giving $\frac{d\omega}{dt}$ we let $dt\rightarrow 0$ for fixed $\epsilon$; then the terms $\frac {dt}{\epsilon^2}$ disappear  and we obtain the simplified ODE formulation with $\omega=\rho, \rho u, \rho v$ successively,\\

$ \frac{d\omega}{dt}(x,y,t,\epsilon)=\frac{1}{\epsilon}[-(\omega(|u|+|v|))(x,y,t,\epsilon)+(\omega u^+)(x-\epsilon,y,t,\epsilon)+  (\omega u^-)(x+\epsilon,$ \begin{equation} y,t,\epsilon)+(\omega v^+)(x,y-\epsilon,t,\epsilon)+(\omega v^-)(x,y+\epsilon,t,\epsilon)]\end{equation}
which is a mere extension of (10, 11) in the $x$ and $y$ directions, to be completed  by 
\begin{equation} u(x,y,t,\epsilon)=\frac{(\rho u)(x,y,t,\epsilon)}{\rho (x,y,t,\epsilon)},  \  and   \   v(x,y,t,\epsilon)=\frac{(\rho v)(x,y,t,\epsilon)}{\rho (x,y,t,\epsilon)}.\end {equation} 
Formulas (43, 44) will be justified by the proof below showing that the system of ODEs has a global solution that gives a weak asymptotic method for the 2-D pressureless fluid system.\\

\textit{Remark on the  numerical calculation of approximate solutions.} We calculate approximate solutions of (43, 44) by using a standard numerical method for ODEs such as a Euler method or a Runge Kutta method (RK4 has given very good results) for fixed $\epsilon>0$. For the space  disretization we use cells of dimension $\epsilon$ in  each axis direction, in which the physical variables are constant. This implies a discretization in time in which $dt$, therefore  $ \frac {dt}{\epsilon^2}$, is not null. In order to permit  calculations in short duration one cannot have  $ \frac {dt}{\epsilon^2}$ too small. Formula (42) describes exactly the transport between cells provided $|u^{\pm} dt|\leq \epsilon$ and $|v^{\pm} dt|\leq \epsilon$. Therefore the terms   $ \frac {dt}{\epsilon^2}$ in (42) can be kept as a  convenient correction from physics to the standard ODE schemes for (43) in order to improve their efficiency in practice, by allowing rather large values of the space step $dt$ ( $|u^{\pm }dt|\leq \epsilon$ and $|v^{\pm} dt|\leq \epsilon$), extending the results obtained with (43) and very small values of $dt$  since (42) describes the physical situation even when $\frac{dt}{\epsilon}$ (the analog of the Courant-Friedrichs-Lewy number $\frac{\Delta t}{\Delta x}$ used in numerical schemes where $\Delta t$ and $\Delta x$ are respectively the time step and the space step)  does not tend to 0. \\


\textit{Extension of the ODEs to 3-D.} In three space dimension with velocity $(u,v,w)$ any cube is surrounded by 6 face neighbors, 12 edge neighbors and  8 vertex neighbors.  Therefore the exact formula for the transport  extending (42) involves remaining matter of density $\omega$ in the cube of center $(x,y,z)$ and edge $\epsilon$, and matter coming from the 26 neighbor cubes instead of the 8 neighbor  squares in (42). The formula can be easily obtained as exposed above in the 2-D case, as in  the appendix in \cite{ColombeauNMPDE}.
 The  formula of the  ODE is a mere extension of the 1-D formula (43),  for  $\omega=\rho,\rho u,\rho v,\rho w$, \\
\\
$ \frac{d\omega}{dt}(x,y,z,t,\epsilon)=\frac{1}{\epsilon}[-\omega (|u|+|v|+|w|)(x,y,z,t,\epsilon)+(\omega u^+)(x-\epsilon,y,z,t, \epsilon)+(\omega u^-)(x+\epsilon,y,z,t, \epsilon)+(\omega v^+)(x,y-\epsilon,z,t, \epsilon)+(\omega v^-)(x,y+\epsilon,z,t, \epsilon)+$\begin{equation}(\omega w^+)(x,y,z-\epsilon,t, \epsilon)+(\omega w^-)(x,y,z+\epsilon,t, \epsilon)].\end{equation} 

For the ODEs (45) it is easy to check that the a priori estimates of section 3 hold without significative modification as well as the global existence-uniqueness result of section 4 and the proof of weak asymptotic method.\\

 We consider initial conditions $\rho_0 \in L^1(\mathbb{T}^3), \rho_0$  positive and  $\vec{ u}_0 \in L^\infty(\mathbb{T}^3)$. We approximate them in the sense of distributions by a family $\{\rho_{0,\epsilon}, u_{0,\epsilon}, v_{0,\epsilon},w_{0,\epsilon}\}_\epsilon$  of continuous functions on $ \mathcal{C}(\mathbb{T}^3)$  with the properties\\
\begin{equation} \rho_{0,\epsilon}(x,y,z)\geq \epsilon  \ \forall x,y,z, \ \ \|\rho_0-\rho_{0,\epsilon}\|_{L^1(\mathbb{T}^3)} \rightarrow 0\end{equation}
  and
\begin{equation} \|u_{0,\epsilon}-u_0\|_{L^\infty(\mathbb{T}^3)}, \|v_{0,\epsilon}-v_0\|_{L^\infty(\mathbb{T}^3)}, \|w_{0,\epsilon}-w_0\|_{L^\infty(\mathbb{T}^3)}  \rightarrow 0\end{equation} when $\epsilon\rightarrow 0$.
From (46) note that $\epsilon$ here plays the role of $m$ in (14), which permits to consider void regions in the initial condition. In the same spirit note that the presence of $\epsilon$ in denominators of (10, 11) permits at the limit $\epsilon\rightarrow 0$ the presence of     concentrations of matter in the solution. The velocity remains bounded and indeed it satisfies the maximum principle in each direction. Summarizing results from the previous sections  we have obtained:\\

\textbf{Theorem 2.}\textit{ Under  the above assumptions  on the initial conditions the system of ODEs (45) complemented with (44) in $u,v,w$ provides a weak asymptotic method of order 1 for  the system  (1, 2) of 3-D pressureless  fluids,  
which is global in space and positive time.}\\

The proof is similar to the 1-D proof and applies also to the pressureless energy equation  (3).\\


\textbf{7. 1-D self-gravitating pressureless fluids.}    The presence of self-gravitation according to Newton's law invalidates the maximum principle in velocity. This difficulty will be easily solved in 1-D thanks to an a priori estimate of the increase in velocity which does not hold in 2-D and 3-D, requesting a more elaborate proof.  From (5-7) the 1-D equations including gravitation are
\begin{equation}\rho_t+(\rho u)_x=0,\end{equation}
\begin{equation}(\rho u)_t+(\rho u^2)_x+\rho \Phi_x=0,\end{equation}
\begin{equation}\Phi_{xx}=4\pi G \rho.\end{equation}
 We assume an  initial condition $\rho_0\in L^1(\mathbb{T})$, with  finite velocity $u_0\in L^\infty(\mathbb{T})$.  \\

We introduce the ODEs   by modifying the Euler equation (11) as follows to take gravitation into account
\begin{equation}\frac{d}{dt}(\rho u)(x,t,\epsilon)=\frac{1}{\epsilon}[(\rho u u^+)(x-\epsilon,t,\epsilon)-(\rho u |u|)(x,t,\epsilon)+(\rho u  u^-)(x+\epsilon,t,\epsilon)]-\rho(x,t,\epsilon)\Phi_x(x,t,\epsilon),\end{equation}
where 

\begin{equation} \Phi_x(x,t,\epsilon)=const+4\pi G\int_{-\pi}^x\rho(\xi,t,\epsilon)d\xi\end{equation}
from (50). The value $const$ in (52) depends on the boundary conditions and does not play any role in the proofs.\\
\\
$\bullet$\textit{Assumptions.}
 We  choose the approximations $\rho_{0, \epsilon} \in \mathcal{C}(\mathbb{T})$ such that \begin{equation} \rho_{0,\epsilon}(x)>0 \ \ \forall x\end{equation}
and \begin{equation}\|\rho_0-\rho_{0,\epsilon}\|_{L^1( \mathbb{T})}\rightarrow 0.\end{equation} 
We assume $u_{0,\epsilon}\in \mathcal{C}(\mathbb{T})$ and
\begin{equation}  \|u_0-u_{0,\epsilon}\|_{L^\infty(\mathbb{T})} \rightarrow 0.\end{equation}

$\bullet$\textit{A priori inequalities}. For fixed $\epsilon$ we assume the existence of a solution  
\begin{equation}[0,\delta(\epsilon)[ \ \ \longmapsto (\mathcal{C}(\mathbb{T}))^2  \  \  \ \ \ \end{equation}
 $$ \  \  \  \ \ \  \  \ \ \ \ \  \ \ \ \ \ \ \ \ \ \ \ \  t \longmapsto [x\mapsto(\rho(x,t,\epsilon),(\rho u)(x,t,\epsilon))]$$
 continuously differentiable such that 
\begin{equation} \exists m>0 \  /  \ \rho(x,t,\epsilon) \geq m \ \forall x\in \mathbb{T}, \end{equation}
\begin{equation} \exists M>0 \ / \ |u(x,t,\epsilon)| \leq M, \  \rho(x,t,\epsilon) \leq M  \ \forall x\in \mathbb{T}. \end{equation}
\\
The proof of (18) gives $\int_{\mathbb{T}}\rho (x,t,\epsilon)dx=\int_{\mathbb{T}}\rho_{0,\epsilon} (x)dx. $ 
 Therefore,  if \begin{equation}K= 8\pi G\int_{\mathbb{T}}\rho_0(x,\epsilon)dx+2.const,\end{equation} from  (52), \begin{equation}|\Phi_x(x,t,\epsilon)|\leq \frac{K}{2} \ \forall (x,t,\epsilon).\end{equation}Let \begin{equation}A(\epsilon):=\|u_{0,\epsilon}\|_\infty+K\delta(\epsilon).\end{equation}  Now we will obtain bounds on $u$ and $\rho$ that depend only on the initial conditions and $\delta(\epsilon)$. Note that in presence of gravitation the velocity can increase with time.\\

\textbf{ Proposition 3.}\textit{  From the above assumptions $\forall x\in \mathbb{T}$ and $\forall t\in [0,\delta(\epsilon)[ $  one has:} 
\begin{equation}  |u(x,t,\epsilon)| \leq A(\epsilon), \ \ \ \ \ \ \ \ \ \ \ \ \ \ \ \ \ \ \ \ \ \ \ \ \ \ \ \ \ \ \ \ \ \ \ \ \ \ \ \ \ \ \ \ \ \ \ \ \ \ \ \ \ \ \ \ \ \ \ \ \ \ \ \ \ \ \ \ \ \ \end{equation}
\begin{equation} \rho_{0,\epsilon}(x) exp(-\frac{A(\epsilon)}{\epsilon}t) \leq \rho(x,t,\epsilon) \leq \|\rho_{0,\epsilon}\|_\infty exp(\frac{2}{\epsilon}const(\epsilon)A(\epsilon)t). \ \ \ \ \ \ \ \ \ \ \ \ \ \end{equation}
\\  
 \textit{Proof}. The proof is an adaptation of the proof of Proposition 1.
From (51) and the mean value theorem\\

 $(\rho u)(x,t+dt,\epsilon)=[\frac{dt}{\epsilon}(\rho u u^+(x-\epsilon,t,\epsilon)+(1-\frac{dt}{\epsilon}|u(x,t,\epsilon)|)(\rho u)(x,t,\epsilon)+$ \begin{equation}\frac{dt}{\epsilon}(\rho u u^-(x+\epsilon,t,\epsilon)]-dt\rho(x,t,\epsilon)\Phi_x(x,t,\epsilon) + dt o_3(x,t,\epsilon,dt)\end{equation}
where $\|o_3(.,t,\epsilon,dt)\|_\infty$ tends to $0$  when $dt$ tends to $0$.
Multiplication of $(\rho u)(x,t+dt,\epsilon)$ and $ \frac{1}{\rho(x,t+dt,\epsilon)}$ from (20) gives $$
\frac{\rho u}{\rho}(x,t+dt,\epsilon)=\frac{\frac{dt}{\epsilon}(\rho u u^+)(x-\epsilon,t,\epsilon)+[1-\frac{dt}{\epsilon}|u|(x,t,\epsilon)](\rho u)(x,t,\epsilon)+\frac{dt}{\epsilon}(\rho u u^-)(x+\epsilon,t,\epsilon)}{\frac{dt}{\epsilon}(\rho u^+)(x-\epsilon,t,\epsilon)+[1-\frac{dt}{\epsilon}|u|(x,t,\epsilon)]\rho(x,t,\epsilon)+\frac{dt}{\epsilon}(\rho u^-)(x+\epsilon,t,\epsilon)}-$$\begin{equation}dt \frac{\rho(x,t,\epsilon)}{\rho(x,t+dt,\epsilon)}\Phi_x(x,t,\epsilon)+dt.o_4(x,t,\epsilon,dt),\end{equation}
\\
where $o_4(x,t,\epsilon,dt) \rightarrow 0$ uniformly in $x$  and $t$ when $dt\rightarrow 0$.  Finally one obtains, as in the passage from (21) to (22), with $K$ defined in (59), 
 \begin{equation} \|u(.,t+dt,\epsilon)\|_\infty\leq  \|u(.,t,\epsilon)\|_\infty+dtK+dt.\|o_4(.,t,\epsilon,dt)\|_\infty\end{equation}
 since $\frac{\rho(x,t,\epsilon)}{\rho(x,t+dt,\epsilon)}<2$ for $dt>0$ small enough depending on $\epsilon$ and $\delta'<\delta (\epsilon)$ if $t\leq\delta'$.\\

Let $f:\mathbb{R}^+\longmapsto\mathbb{R}^+$ be a function such that
\begin{equation} f(t+dt)\leq f(t)+K dt+dt.o(t ,dt)\end{equation}
where $o(t ,dt)\rightarrow 0$ uniformly in $t$ when $dt\rightarrow 0$. Then 
\begin{equation} \forall \tau>0 \ \  f(t+\tau)\leq f(t) +K\tau.\end {equation}
To prove (68) it suffices to apply lemma 1 with $g(t)=f(t)-Kt$.\\

Applying (68) to (66) with $t=0$ and $\tau=t$  for  $t\in [0,\delta']$, then letting $\delta'\rightarrow \delta(\epsilon)$ one has 
\begin{equation}\|u(.,t,\epsilon)\|_\infty\leq \|u_{0,\epsilon}\|_\infty+K\delta(\epsilon)=A(\epsilon)\end{equation}
which replaces here the bound (16). The proof of (63) is identical to the one of Proposition 1. 
For fixed $\epsilon$ one obtains from these estimates a global existence-uniqueness result for system (10, 51, 12, 50) as in section 4.\\

 The proof that the scheme provides a weak asymptotic method of order one in $\epsilon$ is the same as the one in section 3 since the supplementary term $\rho \Phi_x$ from the equation (49) disappears in the proof, by simplification with the corresponding term of the ODE (51). We have obtained the following result.\\

 We consider   initial conditions $\rho_0 \in  L^1(\mathbb{T}),\rho_0(x)\geq 0 \ \forall x $ and  $  u_0 \in L^\infty(\mathbb{T})$. We  approximate them  by a family $\{\rho_{0,\epsilon} \in \mathcal{C}(\mathbb{T}), u_{0,\epsilon}\in \mathcal{C}(\mathbb{T})\}_\epsilon$ which  satisfies  the properties (53-56). Then\\

\textbf{Theorem 3}. \textit{Under these  assumptions  on the initial conditions the system of ODEs (10, 51) complemented by (12, 52) provides a weak asymptotic method of order 1 for  the system  of 1-D self-gravitating pressureless  fluids (48-50)
which is global in space and positive time.}\\


\textbf{8. 2-D and 3-D self-gravitating pressureless fluids.}    In this section we extend the method to 2-D and 3-D. A difficulty stems from the fact that the gradient $\vec{\nabla \Phi}$ of the gravitation potential $\Phi$ is unbounded on point concentrations of matter in 2-D and on point and string concentrations of matter in 3-D, which does not permit to obtain a priori bounds on velocity independent on $\epsilon$ as in (61, 62) in which $\|u_{0,\epsilon}\|_\infty$ can be chosen independent on $\epsilon$ and in which $\delta(\epsilon)$ is a time which is proved independent on $\epsilon$ from the global existence result in section 7. To regularize these singularities we state the Poisson equation by means of a convolution: this can be physically justified from the absence of verification of Newton's law of gravitation at small distances. Then we can obtain a weak asymptotic method in 2-D and 3-D.   Instead of the  1-D integration   (52), we will use the explicit forms of the Green functions  of the Laplace operator in 2-D, i.e. $\Phi(r)=const.log (r), \  r=\sqrt{x^2+y^2}$, and in 3-D, i.e. $\Phi(r)=\frac{const}{r}$.   To shorten the exposition we consider mainly the 2-D case since the 3-D case is  similar.\\ 

The continuity equation is stated as equation (43) with $\omega=\rho$. The Euler equation in the first component $u$ of the velocity vector $(u,v)$
\begin{equation} (\rho u)_t+(\rho u^2)_x +(\rho uv)_y +\rho \Phi_x=0\end{equation}
is formulated as 
$$ \frac{d}{dt}(\rho u)(x,y,t,\epsilon)=\frac{1}{\epsilon}[-((\rho u)(|u|+|v|))(x,y,t,\epsilon)+(\rho u u^+)(x-\epsilon,y,t,\epsilon)+$$  $  (\rho u u^-)(x+\epsilon, y,t,\epsilon)+(\rho u v^+)(x,y-\epsilon,t,\epsilon)+(\rho u v^-)(x,y+\epsilon,t,\epsilon)] - $ \begin{equation}\rho(x,y,t,\epsilon)((\Phi_\epsilon)_x)(x,y,t,\epsilon),\end{equation}
and 
$$ \frac{d}{dt}(\rho v)(x,y,t,\epsilon)=\frac{1}{\epsilon}[-((\rho v)(|u|+|v|))(x,y,t,\epsilon)+(\rho v u^+)(x-\epsilon,y,t,\epsilon)+  (\rho v u^-)$$
$$(x+\epsilon, y,t,\epsilon)+(\rho v v^+)(x,y-\epsilon,t,\epsilon)+(\rho v v^-)(x,y+\epsilon,t,\epsilon)] - \rho(x,y,t,\epsilon)((\Phi_\epsilon)_y)(x,y,t,\epsilon)$$
 for the second Euler equation $(\rho v)_t+(\rho u v)_x +(\rho v^2)_y+\rho \Phi_y=0 $. The velocity $(u,v)$ is given by equations (44). \\

The Poisson equation (7) is stated as follows.
We consider a  regularizing function $\phi\in \mathcal{C}_c^\infty(\mathbb{R}^2) \ (\mathcal{C}_c^\infty(\mathbb{R}^3)$ in the 3-D case),  positive and such that $\int \phi(x,y) dxdy=1$. As usual we set $\phi_{\epsilon^\alpha}(x,y)=\frac{1}{\epsilon^{2\alpha}}\phi(\frac{x}{\epsilon^{\alpha}},\frac{y}{\epsilon^{\alpha}})$ for some $\alpha, \  0<\alpha<1$, to be chosen small enough. We state the Poisson equation in the form
\begin{equation}\Delta(\Phi_\epsilon)=4\pi G .(\rho_\epsilon*\phi_{\epsilon^\alpha})\end{equation}
and we will consider the solution given by the classical Newtonian potentials in 2-D and 3-D.\\


Starting with a $L^1$ initial condition $\rho_0$, the variable $\rho_\epsilon$ is $L^1$ in $x$ with bound independent of $\epsilon$ and $t$, as  proved above in (18).  The notation $const$ will serve for various constant values independent on $\epsilon$.
 Therefore, for given $\phi$,  one has
 \begin{equation} \|\rho_\epsilon*\phi_{\epsilon^\alpha} \|_\infty \leq   \|\rho_\epsilon\|_{L^1}. const. \frac{1}{\epsilon^{2\alpha}} \leq  const. \frac{1}{\epsilon^{2\alpha}}\end{equation} since  $\|\rho_\epsilon\|_{L^1}$ is bounded independently of $\epsilon$. The factor $\frac{1}{\epsilon^{2\alpha}}$ comes from the formula $\phi_{\epsilon^\alpha}(x,y)=\frac{1}{\epsilon^{2\alpha}}\phi(\frac{x}{\epsilon^\alpha},\frac{y}{\epsilon^\alpha})$ as a mollifier in 2-D. \\
 
  From the  2-D Green function of the Laplace operator one obtains the classical formula of the Newtonian potential $\Phi(x,y)$, in which $\rho$ denotes any continuous density of matter on $\mathbb{T}^2$ (that will be  applied later to $\rho_\epsilon*\phi_{\epsilon^\alpha}$ in place of $\rho$): 
 
 $$\Phi(x,y)=const.\int_{\mathbb{T}^2} \rho(\xi,\eta)  log([(x-\xi)^2+(y-\eta)^2])d\xi d\eta.$$
 \\
 By differentiation
 
 $$\Phi_x(x,y)=const.\int_{\mathbb{T}^2} \rho(\xi,\eta)  \frac{x-\xi}{(x-\xi)^2+(y-\eta)^2}d\xi d\eta,$$
 which implies 
 $$|\Phi_x(x,y)|+|\Phi_y(x,y)|\leq const.\int_{\mathbb{T}^2} \rho(\xi,\eta)  \frac{1}{\sqrt{(x-\xi)^2+(y-\eta)^2}}d\xi d\eta.$$
We split this integral into 
 $$|\Phi_x(x,y)|+|\Phi_y(x,y)|\leq const. \int_{\mathbb{T}^2,r>1} \rho(\xi,\eta) d\xi d\eta+const.\int_{\mathbb{T}^2, r<1} \|\rho\|_\infty\frac{1}{r} r dr d\theta,$$ 
with the notation $r=\sqrt{(x-\xi)^2+(y-\eta)^2}$,  i.e.
$$|\Phi_x(x,y)|+|\Phi_y(x,y)|\leq const. \| \rho\|_{L^1(\mathbb{T}^2)}  +const. \|\rho\|_{L^\infty(\mathbb{T}^2)}.$$  
Now we apply this bound in which $\rho$ is replaced by $\rho_\epsilon*\phi_{\epsilon^\alpha}$. From the convolution $\|\rho_\epsilon*\phi_{\epsilon^\alpha}\|_{L^1(\mathbb{T}^2)}\leq\|\rho_\epsilon\|_{L^1(\mathbb{T}^2)}\leq const$ and $\|\rho_\epsilon*\phi_{\epsilon^\alpha}\|_\infty\leq\frac{const}{\epsilon^{2\alpha}}$ from (73). We finally obtain 

\begin{equation}|\Phi_x(x,y)|+|\Phi_y(x,y)|\leq const.\frac{1}{\epsilon^{2\alpha}}. \end{equation}
The bound (74) will replace the bound (60) $\|\Phi_x(.,t,\epsilon)\|_\infty \leq const \ \forall\epsilon$ of the 1-D case. One obtains an analog of (62)  but with a factor $\frac{1}{\epsilon^{2\alpha}}$ from (74) i.e. $|u(x,t,\epsilon)|\leq  \|u_{0,\epsilon}\|_\infty+\frac{const.\delta(\epsilon)}{\epsilon^\alpha} $. The bound $K$ in (66) and (69) is  replaced by $ \frac{const}{\epsilon^{2\alpha}}$  (this replacement holds also in the sequel).  The variable $\rho u$  is no longer $L^1$ uniformly in $\epsilon$ as in the cases considered up to now in which the velocity was bounded. We have  $\|\rho u\|_{L^1(\mathbb{T}^2)}\leq \frac{const}{\epsilon^{2\alpha}}$.\\

 This  allows the possibility of  infinite velocity, as exposed in \cite{Saari} concerning the 3-D N-body problem, which makes a great difference with the 1-D case,  at the same time as this does not significantly perturb the proof of weak asymptotic method since a bound $O(\epsilon)$ previously obtained in the 1-D case for the continuity and Euler equations will be simply replaced by a bound $O(\epsilon^{1-2\alpha})$ and one can choose $\alpha>0$ as small as needed.\\


 The 1-D proofs of weak asymptotic method for the continuity equation (5) and the Euler equation (7) apply in 2-D and 3-D   for $\alpha>0$ small enough since the above bound $\frac{const}{\epsilon^{2\alpha}}$ (respectively $\frac{const}{\epsilon^{3\alpha}}$) is in factor of $\epsilon$ and one can choose $\alpha>0$ as small as needed. For the Poisson equation (8) one has to prove that $$ \forall\psi\in \mathcal{C}_c^\infty (\mathbb{R}^3) \ \int(\Delta\Phi_\epsilon(x,y,z,t)-4\pi G\rho(x,y,z,t,\epsilon)) \psi(x,y,z) dxdydz\rightarrow 0$$ when $\epsilon\rightarrow 0$ which follows from (72) and the $L^1$ stability in $\rho$. Therefore here in 2-D the weak asymptotic method is of order  $1-2\alpha$ for the continuity and Euler equations and $2\alpha$ for the Poisson equation.\\

The proof extends to 3-D from a similar statement of the Poisson equation and a similar analysis of the Newtonian potential: one obtains \\

$|\Phi_x(x,y,z)|+|\Phi_y(x,y,z)|+|\Phi_z(x,y,z)| \leq $\\

 $const.\int_{\mathbb{R}^3}\rho(\xi,\eta,\mu)\frac{1}{(x-\xi)^2+(y-\eta)^2+(z-\mu)^2}d\xi d\eta d\mu 
\leq $ \\ 

$ const.\int_{r>1}\rho(\xi,\eta,\mu)d\xi d\eta d\mu +const.\int_{r<1}\|\rho\|_\infty dr  \leq const.(\|\rho\|_{L^1}+\|\rho\|_\infty)$. \\

In 3-D one obtains
a weak asymptotic method of order $1-3\alpha$ for the constitutive and Euler equations, choosing $\alpha<\frac{1}{3}$, and $3\alpha$ for the Poisson equation.\\

 Finally we obtain: let  be given initial conditions $\rho_0 \in L^1(\mathbb{T}^3), \ \ \rho_0$  positive with possibility of void regions and $ \ u_0 \in L^\infty(\mathbb{T}^3)$. We assume that the approximations   $\{\rho_{0,\epsilon}, u_{0,\epsilon}\}$ satisfy $\rho_{0,\epsilon}(x)>0 \ \ \forall x\in \mathbb{T}^3, \ \|\rho_0-\rho_{0,\epsilon}\|_{L^1(\mathbb{T}^3)}\rightarrow 0$  and $  \|u_0-u_{0,\epsilon}\|_{L^\infty(\mathbb{T}^3)}\rightarrow 0.$\\

\textbf{Theorem 4.}\textit{ Under  the above assumptions  on the initial conditions the system of ODEs (45 in $\rho$, 71 in $\rho\vec{U}$) complemented with (44, 72) provides a weak asymptotic method  for  the system  (5-7) of 2-D and 3-D self-gravitating pressureless  fluids
which is global in space and positive time.}\\
\\

\textbf{9. The case of a finite total mass on $\mathbb{R}^n.$} In this section we sketch how the constructions can be modified when the torus $\mathbb{T}^n$ is replaced by the euclidean space $\mathbb{R}^n$ on which the initial condition is a finite mass. We denote by  $\mathcal{C}_b(\mathbb{R})$ the Banach space of continuous bounded functions on $\mathbb{R}$.\\
 
$\bullet$\textit{Assumptions.} We assume that
$\rho_0\in L^1(\mathbb{R}) $
 and that we can choose the approximations $\rho_{0, \epsilon} \in L^1(\mathbb{R}) \cap \mathcal{C}_b(\mathbb{R})$  such that $\forall \epsilon>0$ small enough there is an auxiliary  function
 
 \begin{equation} p_{\epsilon}\in L^1(\mathbb{R})\cap \mathcal{C}_b(\mathbb{R}), \ \ p_{\epsilon}(x)>0 \  \ \forall x\end{equation}
 such that 
\begin{equation}\exists a(\epsilon),b(\epsilon)>0 \ /  \ a(\epsilon) p_{\epsilon}(x)<\rho_{0,\epsilon}(x) <b(\epsilon) p_{\epsilon}(x) \ \ \forall x \in \mathbb{R}\end{equation}
 \\
and \begin{equation}\|\rho_0-\rho_{0,\epsilon}\|_{L^1( \mathbb{R})}\rightarrow 0.\end{equation} 

The assumption on $u_0$ and its approximations is  
\begin{equation} \ u_0\in L^\infty(\mathbb{R}) \ \   and \ \  \|u_0-u_{0,\epsilon}\|_\infty \rightarrow 0.\end{equation}
Further, we assume that $p_{\epsilon}$ satisfies 
\begin{equation} \frac{p_{\epsilon}(x\pm\epsilon)}{p_{\epsilon}(x)} \leq const(\epsilon) \ \forall x\in \mathbb{R}.\end{equation}
This last assumption is satisfied if $p_{\epsilon}(x)=exp(-a|x-x_0|)$, for  $a>0$  and any $x_0\in \mathbb{R}$ (in this case $p_{\epsilon}$ does not depend on $\epsilon$) but not if $p_{\epsilon}(x)=exp(-x^2)$. It is immediate that the assumptions (75-77, 79) can be satisfied if $\rho_0\in L^1(\mathbb{R})\cap L^\infty(\mathbb{R})$ has compact support. \\
\\
$\bullet$\textit{The Banach space.} We consider the Banach space
 $$E_{p_{\epsilon}}=\{X\in \mathcal{C}_b(\mathbb{R}) / \exists \lambda>0 / |X(x)| \leq \lambda p_{\epsilon}(x) \forall x \in \mathbb{R}\}$$ equipped with the norm 
\begin{equation} \|X\|_{p_{\epsilon}}=sup_{x\in \mathbb{R}}\frac{|X(x)|}{p_{\epsilon}(x)}=inf\{\lambda>0 / |X(x)|\leq \lambda p_{\epsilon}(x) \forall x \in \mathbb{R}\}.\end{equation}
$\bullet$\textit{A priori inequalities}. For fixed $\epsilon>0$ we assume the existence of a solution  
$$[0,\delta(\epsilon)[ \ \ \longmapsto (E_{p_{\epsilon}})^2  \  \  \ \ \ $$
 $$ \  \  \  \ \ \  \  \ \ \ \ \  \ \ \ \ \ \ \ \ \ \ \ \  t \longmapsto [x\mapsto(\rho(x,t,\epsilon),(\rho u)(x,t,\epsilon))]$$
 continuously differentiable such that 
\begin{equation} \exists m>0 \  /  \ \rho(x,t,\epsilon) \geq m p_{\epsilon}(x) \ \forall x\in \mathbb{R}, \end{equation}
\begin{equation} \exists M>0 \ / \ |u(x,t,\epsilon)| \leq M, \  \rho(x,t,\epsilon) \leq M p_{\epsilon}(x) \ \forall x\in \mathbb{R}. \end{equation}
 \\
We sketch the proof in the case of the system (48-50) of 1-D selfgravitating fluids. One obtains, with $A(\epsilon)$ defined in (61),\\

\textbf{ Proposition 4.}\textit{  $\forall x\in \mathbb{R},\forall t\in [0,\delta(\epsilon)[$  one has} 
\begin{equation}  |u(x,t,\epsilon)| \leq A(\epsilon) \ \ \ \ \ \ \ \ \  \ \ \ \ \ \ \ \ \ \ \ \ \ \ \ \ \ \ \ \ \ \ \ \ \ \ \    \ \ \ \ \ \ \ \ \ \ \ \ \ \  \ \ \ \ \ \ \ \ \ \ \ \ \ \  \end{equation} and
\begin{equation} a(\epsilon)p_{\epsilon}(x) exp(-\frac{A(\epsilon)}{\epsilon}t) \leq \rho(x,t,\epsilon) \leq b(\epsilon) p_{\epsilon}(x) exp(\frac{2}{\epsilon}const(\epsilon)A(\epsilon)t). \end{equation}
 
 \textit{Proof}. The proof is an adaptation of the proof of Proposition 3 given on $\mathbb{T}$  therefore we only sketch the changes. Equation (10) implies again (19) from the mean value theorem in the Banach space $E_{p_\epsilon}$ but here $\|o(.,t,\epsilon,dt)\|_{p_\epsilon} \rightarrow  0$ instead of the sup norm $\|o(.,t,\epsilon,dt)\|_\infty $ in (19). The inversion gives \\
 
 $ \frac{1}{\rho(x,t+dt,\epsilon)}=$\\
 $[\frac{dt}{\epsilon}(\rho u^+)(x-\epsilon,t,\epsilon)+[1-\frac{dt}{\epsilon}|u|(x,t,\epsilon)]\rho(x,t,\epsilon)+\frac{dt}{\epsilon}(\rho u^-)(x+\epsilon,t,\epsilon)]^{-1}+dt.o_1(x,t,\epsilon,dt)$\\
 \\
where $|o_1(x,t,\epsilon,dt)|\leq const.\frac {o_2(x,t,\epsilon,dt)}{p_{\epsilon}(x)}$ for some $o_2$. This follows from  the formula $\frac{1}{a+\mu}=\frac{1}{a}-\frac{\mu}{a^2}+o(\frac{\mu}{a^2})$ for $a>0$ and small $\mu$. Here $a$ is the sum of the first three terms in the second member of (19). Therefore, from (81), $\exists \alpha>0$ such that $a\geq\alpha p_{\epsilon}(x)$.
In (19) one has  $\mu=dt.o(x,t,\epsilon,dt)$. Therefore $|\mu(x,t,\epsilon,dt)|\leq \|o(.,t,\epsilon,dt)\|_{p_\epsilon}.p_{\epsilon}(x)$ from (80). 
Then $|o_1(x,t,\epsilon,dt)|=|\frac{-\mu}{a^2}+o(\frac{\mu}{a^2})|\leq \frac{\|o_\epsilon\|_{p_\epsilon}p_\epsilon(x)}
{(\alpha p_\epsilon(x))^2} $.\\

From (49) and the mean value theorem in the Banach space $E_{p_\epsilon}$\\
\  \\
 $(\rho u)(x,t+dt,\epsilon)=[\frac{dt}{\epsilon}(\rho u u^+(x-\epsilon,t,\epsilon)+(1-\frac{dt}{\epsilon}|u(x,t,\epsilon)|)(\rho u)(x,t,\epsilon)+\frac{dt}{\epsilon}(\rho u u^-(x+\epsilon,t,\epsilon)]-dt\rho(x,t,\epsilon)\Phi_x(x,t,\epsilon) + dt o_3(x,t,\epsilon,dt)$\\
\ \\
where $\|o_3(.,t,\epsilon,dt)\|\rightarrow 0$  when $dt$ tends to $0$.\\

Multiplication of $(\rho u)(x,t+dt,\epsilon)$ and $ \frac{1}{\rho(x,t+dt,\epsilon)}$ makes the $p_{\epsilon}(x)$  disappear and  one obtains (66) with the sup norm for the remainder  in velocity.\\

 The upper bound of $\|\rho(.,t,\epsilon)\|_{p_{\epsilon}}$ is obtained as follows. From the definition (80) of the norm $\| .\|_{p_{\epsilon}}$ and  from the property (79) of $p_{\epsilon}$, one has successively 
 $$ \rho(x\pm \epsilon,t,\epsilon)\leq \|\rho(.,t,\epsilon)\|_{p_{\epsilon}}p_{\epsilon}(x\pm\epsilon).$$
 and
  $$ \rho(x\pm \epsilon,t,\epsilon)\leq \|\rho(.,t,\epsilon)\|_{p_{\epsilon}}const(\epsilon)p_{\epsilon}(x).$$
  The ODEs (10) and (62) give $$\frac{d\rho(x,t,\epsilon)}{dt}\leq \frac{1}{\epsilon}[\rho(x-\epsilon,t,\epsilon)+\rho(x+\epsilon,t,\epsilon)]A(\epsilon).$$ Therefore 
  $$\frac{d\rho(x,t,\epsilon)}{dt}\leq \frac{2}{\epsilon}\|\rho(.,t,\epsilon)\|_{p_{\epsilon}}.const(\epsilon)p_{\epsilon}(x)A(\epsilon).$$
  By integration 
  $$\rho(x,t,\epsilon)\leq\rho_0(x,\epsilon)+\frac{2}{\epsilon}const(\epsilon)p_{\epsilon}(x)A(\epsilon)\int_0^t\|\rho(.,s,\epsilon)\|_{p_{\epsilon}}ds.$$
  Using (80) and dividing by $p_{\epsilon}(x)>0$ one obtains 
  $$\frac{\rho(x,t,\epsilon)}{p_{\epsilon}(x)}\leq \|\rho_{0,\epsilon}\|_{p_{\epsilon}} +\frac{2}{\epsilon}const(\epsilon)A(\epsilon)\int_0^t\|\rho(.,s,\epsilon)\|_{p_{\epsilon}}ds \ \forall x\in \mathbb{R}.$$
  Finally, using (80) again,

  $$\|\rho(.,\epsilon,t)\|_{p_{\epsilon}}\leq \|\rho_{0,\epsilon}\|_{p_{\epsilon}} +\frac{2}{\epsilon}const(\epsilon)A(\epsilon)\int_0^t\|\rho(.,s,\epsilon)\|_{p_{\epsilon}}ds $$
  we conclude by applying  Gronwall formula. \\

\textbf{10. Connection with the N-body problem}. The system of PDEs (5-7) is a continuous extension of the discrete N-body problem  by stating the two Newton's laws in continuous form: Newton's law of motion gives the Euler equation and  Newton's law of gravitation gives the Poisson equation. Not only it has been noticed that the limit of the weak asymptotic method can give concentrations of matter but also one can notice  that all theorems and proofs in this paper hold without change if the initial condition $\rho_0$ is a positive bounded  Radon measure in place  of a $L^1$  function, provided the velocity is not discontinuous on the concentration points of $\rho_0$ (if not one has to be cautious see \cite{ColombeauSiam} p. 1911), since a bounded Radon measure $\rho_0$ can be approximated by continuous functions $\rho_{0,\epsilon}$ with properties (53-55), which permits to consider N 
point-bodies as initial condition in the weak asymptotic method. \\
\\
Then setting 
\begin{equation}\rho(x,y,z,t,\epsilon)=\sum_i m_i\delta(x-X_i(t,\epsilon),y-Y_i(t,\epsilon),z-Z_i(t,\epsilon),t,\epsilon),\end{equation}
\begin{equation}(\rho\vec{U})(x,y,z,t,\epsilon)=\sum_i m_i \vec{U}_i(t,\epsilon)\delta(x-X_i(t,\epsilon),y-Y_i(t,\epsilon),z-Z_i(t,\epsilon),t,\epsilon),\end{equation}
where $\delta(M,t,\epsilon)$ a regularization by means of the parameter $\epsilon$ of the 3-D Dirac delta Radon measure. Inserting (85, 86) into the system of ODEs (45, 71) one obtains formally at the limit $\epsilon=0$ the two respective classical equations 
\begin{equation} \frac{d}{dt} \vec{r}_i=\vec{U}_i,  \  \frac{d}{dt} \vec{U}_i=-G\sum_{j\not=i} m_j\frac{\vec{r}_i-\vec{r}_j}{\|\vec{r_i}-\vec{r_j}\|^3},\end{equation}
using the classical formula of the Newtonian potential $\Phi(\vec{r})=G\int\frac{\rho(\vec{x})}{\|\vec{x}-\vec{r}\|} d^3x$. It appears that the convolution in (72) replaces the classical regularization $|\|\vec{r_i}-\vec{r_j}\|^2+\epsilon^2|^{\frac{3}{2}}$ of the denominators in (87).\\ 

From these remarks  the systems of ODEs in infinite dimension (45 in $\rho$, 71 in $\rho \vec{U}$) completed by the regularized Poisson equation (72), originally issued from system (5-7), can also be viewed as an extension of the classical system (87) of ODEs  of the N-body problem to the continuous case at the limit $\epsilon\rightarrow 0^+$. \\

\vskip 3 cm
\textit{}\\

\textit{ }\\
\\
\\
\\
\\
\\
\\
\\
\\
\\
\\
\\
\\
\\
\\
\\
\\
\\
\\
\\
\\
\\
\\
\\
\\
\\
\\
\\
\\
\\
\\

\begin{figure}[h]
\centering
\includepdf[width=\textwidth]{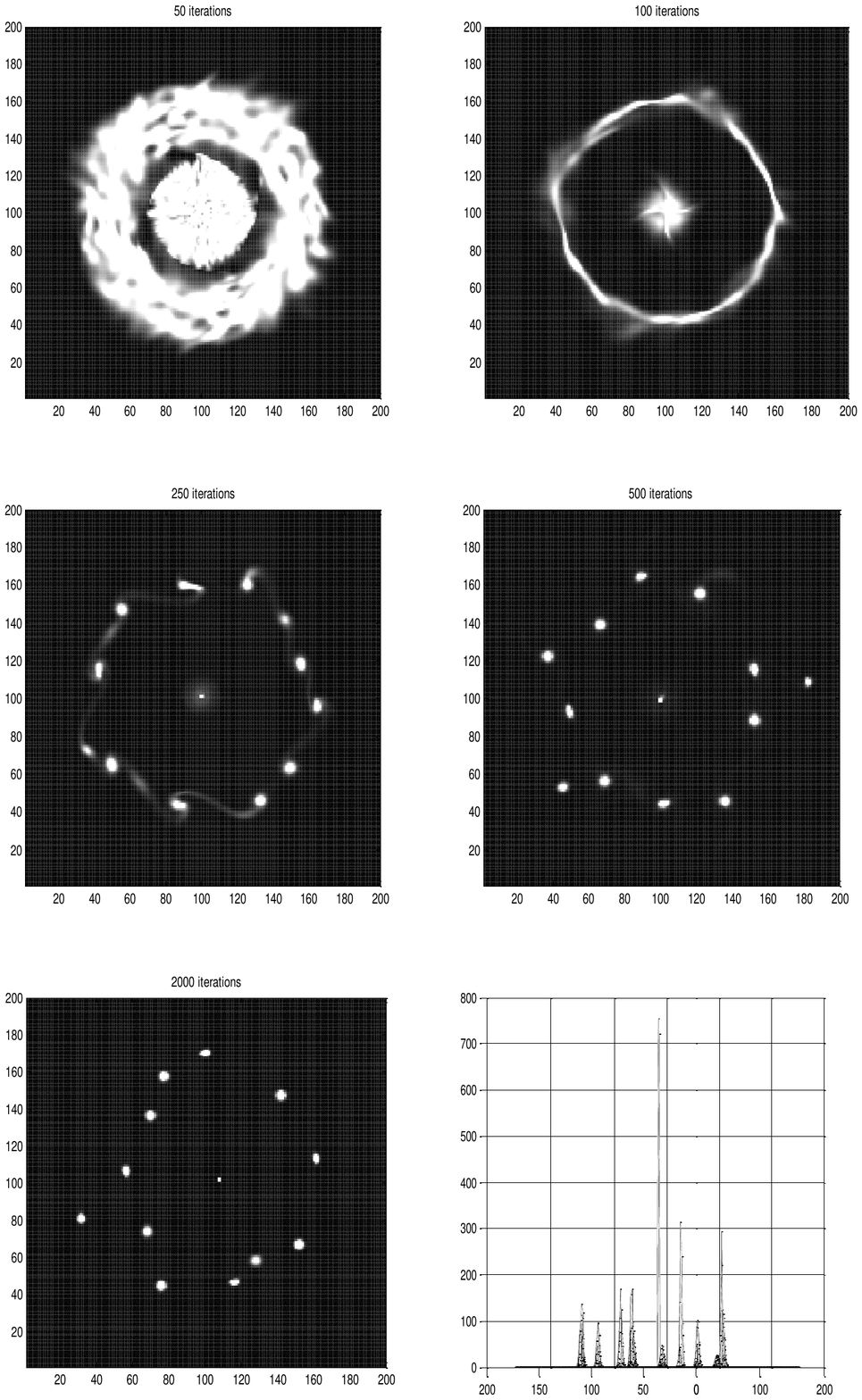}
\end{figure}

Note that the modelling of the $N$-body problem from the weak asymptotic method in this paper is far more refined than the usual one from (87) since the bodies have a volume: consider a solitary concentration of matter such as the "`sun"' in the simulation of section 11 in absence of all planets: a motionless mass located on $[-\epsilon,+\epsilon]^2 \ \forall \epsilon$ is not a singular point of the vector field because of the last term in  (71). Such a solitary mass goes on collapsing towards a point endlessly: this is in agreement with physics in which this gravitational collapse is only stopped by pressure not considered in the equations of this paper, see \cite{Colombeaupressure} for that:  figure 7 in \cite{ColombeauNMPDE} describes the gravitational collapse of a cloud of gas of dimension $>$ Jeans' length \cite{Coles,Peacock} to an equilibrium in which pressure compensates exactly gravitation (a rough simulation of star formation). Considering bodies with a nonzero volume has permitted more evolved collision tests in which one has observed ejections of clouds of matter and of pieces of broken bodies. In this viewpoint the method in this paper is a refinement of  the classical $N$-body problem (87) at the same time it provides a solution process to  the classical system (5-7) of continuous fluid mechanics.\\

\textbf{11. An example of numerical simulation: formation and evolution of a planetary system.} The properties of the ODEs permit to prove convergence of the classical explicit Euler order 1 method. The space discretization is done by cubes of side $\epsilon$  parallel to the axis. When the ODEs (10, 11) are treated with the standard Runge Kutta RK4 method one observes already in 1-D an independence of the numerical result of the value $\frac{dt}{\epsilon}$ where $dt$ denotes the discretization step. Indeed one observes that one obtains in general exactly the same numerical result with very small values of $\frac{dt}{\epsilon}$, such as 
$10^{-4}$- as it should be theoretically since (for the solution of the ODEs) at  first $\epsilon$ is fixed small and $dt$ tends to 0, and only then (for the weak asymptotic method) $\epsilon$ tends to 0 - 
and with  rather  large values of $\frac{dt}{\epsilon}$,  provided they satisfy $\|u\|_\infty\frac{dt}{\epsilon}<1$. In contrast, the scheme in \cite{ColombeauSiam,ColombeauNMPDE} from the explicit Euler order one method for (10, 11) gives a numerical result whose quality diminishes for too small values of the Courant-Friedrichs-Lewy number $r:=\frac{\Delta t}{\Delta x}$. This observation of  the numerical solution of the ODEs (10,11) with the RK4 method is a great advantage:  when some region of calculation imposes a small value of $\frac{dt}{\epsilon}$ (due to a large velocity there) the quality of the whole numerical result is not affected in other regions where the velocity is small.\\

In 2-D one uses the improved discretization (42) involving terms $\frac{dt}{\epsilon^2}$ for a sharper  discretization than that of the ODE (43). In presence of gravitation (system (5-7)) one presents a simulation of the formation of a "planetary system" from a rotating cloud of dust presenting some small heterogeneities (see  realistic simulations of physics in \cite{Stark}). The initial conditions are those in figure 4 of \cite{ColombeauNMPDE}:  initial values of density at random in the initial cloud in form of a disk, initial velocities tangential to circles centered in the center of the window, also with random values, with null velocity in a neighborhood of the center.  One can apply   theorem 4  section 8 which states that the solutions of the ODEs (45, 71) complemented by  (44, 72) tend to satisfy the equations, then solve numerically these ODEs with the explicit Euler order one method. From the top-left panel to the bottom-left panel: the disk of dust separates into a central concentration  that is collapsing to the center and a second concentration into a ring (top-left panel). The central concentration  collapses to a point object, "`the star"' and  the ring becomes thin (top-right panel). Then  "planet formation" in the ring starts at once (middle-left panel). Planets  become distinct objects that under influence of gravitation can slowly change location and take different distances from the star. One observes  rotation with changes of aspect of the global system. If one excludes a few observations of ejection of a planet from the window, the whole system looks rather stable during a reasonable test of a few hours stopped for convenience. One does not know if one observes nonclosed orbits or long periodic orbits. Satellites of  planets cannot be observed because the discretization around the planets is too coarse but their embryos can be perceived in the middle left panel. While the star is concentrated on one cell, the planets are smeared over a few cells due to their weaker mass and some  rotation around their center can be guessed in  middle-left panel  before their complete formation. Since there is some small random choice in the initial conditions so as to produce collapse in the unstable situation of symmetry of rotation of the initial cloud the  observed results can be significantly different. The bottom-right panel shows that in this simulation more than 50 per cent of the matter is concentrated in the star (often up to 80 per cent).  This simulation requires only a few minutes on a standard PC; it can be continued  to observe the slow evolution of the planetary system.\\

 The original mathematical novelty is that now one knows from theorem 4 that the observed numerical results represent an  approximate solution  of  the system  (5-7), with however a lack of a uniqueness result of admissible weak asymptotic methods, which, up to now, did not cause problem  as if some uniqueness did hold true in practice in case of classical tests from computational fluid dynamics such as those in \cite{ColombeauSiam,ColombeauNMPDE, Colombeauideal}. This simulation requires only a few minutes on a standard PC.\\

In the same way one could study a galaxy instead of a planetary system, as an illustration of theorem 4. Then one needs to extend the Euler-Poisson system (5-7) to the case of two species of matter (dark matter and baryonic matter) in expanding background. The Euler-Poisson system with two species of matter ($i=1,2$ in (88, 89)) is simply the juxtaposition of the continuity (88) and Euler (89)  equations of each matter that are linked by the Poisson equation (90) which involves both matters  \cite{Coles} p. 242:
\begin{equation}\frac{\partial \rho_i}{\partial t}+\vec{\nabla}.(\rho_i \vec{u_i})=0, \ \ i=1,2\end{equation} 
\begin{equation}\frac{\partial}{\partial t}(\rho_i \vec{u_i})+\vec{\nabla}.(\rho_i \vec{u_i} \otimes \vec{u_i})
+\rho_i \vec{\nabla \Phi}=\vec{0}, \  \  i=1,2\end{equation}
\begin {equation}\Delta \Phi=4\pi G (\rho_1+\rho_2). \end{equation}
The proofs in sections (8, 9) and theorems 3 and  4 extend at once. Expanding background is treated by the introduction of the  scale factor $a(t)$ that describes the (known from the Friedman equations \cite{Coles} p. 294, \cite{Peacock} p. 463, \cite{Peter} p. 233) expansion of the background and transforms the Euler-Poisson system into a mathematically equivalent system, \cite {Coles} p. 94. Figure 5 in \cite{ColombeauNMPDE} shows a 1-D numerical simulation of system (88-90) from the explicit Euler order one method for the numerical solution of (10, 51, 12, 52) (the 1-D scheme in \cite{ColombeauNMPDE}): at the time of decoupling of radiation and baryonic matter the randomly distributed baryonic matter falls into the potential wells of dark matter, forming the future galaxies.\\

The equations of physics such as those stated in introduction  are ideal equations marred by an uncertainty (due to idealizations, in particular due to the fact that the molecular structure of matter is not taken into account) that can be naturally considered in the sense of distributions in the space variables. The approximate solutions produced by an asymptotic method satisfy the equations modulo this uncertainty for $\epsilon>0$ small enough. In this sense the weak asymptotic methods could  produce physically acceptable results although they do not provide an exact mathematical solution of the ideal equations. \\

 Unfortunately the problem of uniqueness of the limit $\epsilon\rightarrow 0$ from weak asymptotic methods remains unsolved. Uniqueness can be understood in two domains. In the domain of classical tests used to validate numerical schemes it appears that some uniqueness certainly exists. Outside this domain, as noticed in  \cite{ColombeauSiam} p. 1911 and in \cite{Danilov0} p. 7 when the initial condition for density contains the Dirac delta measure,  the classical form of the equation and/or the initial conditions can cover different real processes, so uniqueness should rely on some more precise statements of the equations and initial conditions at some infinitesimal level, as noticed in   \cite{ColombeauSiam} p. 1911 \cite{Danilov0} p. 7:  in the case of an interplay between a concentration in density and a discontinuity in velocity the repartition of velocity inside the conscentration of density governs the aspect of the solution.\\

\textbf{ 12. Conclusion.} We have constructed  with full mathematical proofs  approximate solutions, more precisely weak asymptotic methods, for the general Cauchy problem for the  system of pressureless fluid dynamics in  3-D, possibly in presence of self-gravitation, for which there was   no known mathematical solution  in physically relevant situations, besides the recognized importance of this system. The construction of approximate solutions has been extended to presence of pressure in \cite{Colombeaupressure}.\\

 This has been done by means of  two nonlinear ODEs in Banach spaces (one for the continuity equation, one for the Euler equation) for which we prove existence-uniqueness of global solutions for positive time. Further, to prove the pertinence of our method, we have checked that  numerical solutions of these equations  have always given back in the classical tests the known  exact solutions (or the widely accepted solutions in absence of  exact solutions) see \cite{ColombeauSiam,ColombeauNMPDE,Colombeauideal} for the numerical scheme that inspired the more elaborate method in this paper, see \cite{Colombeaupressure} in presence of pressure. To some extent, these verifications replace the lack of a uniqueness proof of  the limits of the weak asymptotic methods we construct. Indeed the result of these numerical verifications could be expected since the theoretical method in this paper has been obtained as an abstract version of the numerical method in \cite{ColombeauSiam,ColombeauNMPDE,Colombeauideal}  in which many numerical tests are reported.\\

 The numerical simulation of formation and evolution of a planetary system from a rotating disk of dust shows that, from the theorems proved in this paper,  complex  physical phenomena  which are presently of great scientific importance and are widely reproduced today from heuristic numerical simulations by physicists and engineers in computational fluid dynamics can now be attained  by mathematical rigor even in absence of known   mathematical solutions  that would be classical functions or distributions.\\


\begin{thebibliography}{< >}

\bibitem{Albeverio} S. Albeverio, O.S. Rozanova, V.M. Shelkovich. Transport and concentration processes in the multidimensional zero-pressure gas dynamics model with the energy conservation law. arXiv: 1101.581v1 [math-ph] 30 Jan 2011.

\bibitem{Blelloch} G. Blelloch, G. Narlikar. A practical comparison of N-body algorithms. In "`Parallel Algorithms"', DIMACS Series in Discrete Math. and Computer science. American Mathematical Society, Vol. 30, 1997, pp.81-96.
\bibitem{Bouchut} F. Bouchut, S. Jin, X. Li. Numerical approximations of pressureless and isothermal gas dynamics. SIAM J. Numer. Anal. 41, 2003, pp. 135-158.
\bibitem{Charru} F. Charru. Hydrodynamics Instabilities. Cambridge texts in applied mathematics. Cambridge University Press. 2011.
\bibitem{Chertock} A. Chertock, A. Kurganov, Y. Rykov. A new sticky particle method for pressureless gas dynamics. SIAM J. Numer. Anal. 45, 2007, pp. 2408-2441.
\bibitem{Coles} P. Coles, F. Lucchin. Cosmology. The Origin and Evolution of Cosmic Structure. 2002. Wiley, second edition.





\bibitem{ColombeauSiam} M. Colombeau. A method of projection of delta waves in a Godunov scheme and application to pressureless fluid dynamics. SIAM J. Numer. Anal. 48, 5, 2010, pp. 1900-1919.

\bibitem{ColombeauNMPDE} M. Colombeau. A consistent numerical scheme for self-gravitating fluid dynamics. Num. Methods for PDEs. 29, 1, 2013, pp. 79-101.


\bibitem{Colombeauideal} M. Colombeau. A simple  numerical scheme for the 3-D system of ideal gases and a study of approximation  in the sense of distributions. J. Comput. Appli. Math. 248, 2013, pp.15-30.
   
\bibitem{Colombeaupressure} M. Colombeau. Weak asymptotic methods for some systems of fluid dynamics with pressure terms. preprint.


\bibitem{Danilov0} V.G. Danilov. Remarks on vacuum state and uniqueness of concentration process. Electronic J. of Differential  Eqs, 34, 2008, pp.1-10.

\bibitem{Danilov1}  V. G. Danilov, G.A. Omel'yanov, and V.M. Shelkovich. Weak Asymptotic Method and Interaction of Nonlinear Waves, AMS Translations vol 208,  2003, pp 33-164.



\bibitem{Mitrovic} V. G. Danilov, D. Mitrovic. Delta shock wave formation in the case of triangular hyperbolic system of conservation laws. J. Differential Equations 245, 2008, pp. 3704-3734.


\bibitem{Shelkovich2} V. G. Danilov, V.M. Shelkovich. Dynamics of propagation and interaction of $ \delta$  shock waves in conservation law systems. J. Differential Equations 211, 2005, pp. 333-381.

\bibitem{Shelkovich3}  V. G. Danilov, V.M. Shelkovich. Delta-shock wave type solution of hyperbolic systems of conservation laws. Quart. Appl. Math. 63, 2005, pp. 401-427.



\bibitem{Rykov} Weinan  E, Yu.G. Rykov, Ya.G. Sinai. Generalized variational principles, global weak solutions and behavior with random initial data for systems of conservation laws arising in adhesion particle dynamics. Comm. Math. Phys. 177, 1996, pp. 349-380.
\bibitem{Heggie} D. Heggie, P. Hut. The gravitational million-body problem. Cambridge University Press, Cambridge, 2003.


\bibitem{LeVeque} R. J. LeVeque. The dynamics of pressureless dust clouds and delta waves. J. Hyperbolic Diff. Eq. 1, 2004, pp. 315-327.


\bibitem{Nguyen} T. Nguyen, A. Tudorascu. Pressureless Euler/Euler-Poisson systems via adhesion dynamics and scalar conservation laws. SIAM J. Math. Ana. 40, 2, 2008, pp. 754-775.



\bibitem{Nilsson1} B. Nilsson, V.M. Shelkovich. Mass, momentum and energy conservation laws in zero-pressure gas dynamics and $\delta$-shocks. Applicable Analysis, 90,1, 2011, pp. 1677-1689.


\bibitem{Nilsson2} B. Nilsson, O.S. Rozanova, V.M. Shelkovich. Mass, momentum and energy conservation laws in zero-pressure gas dynamics and $\delta$-shocks II. Applicable Analysis, 90,5, 2011, pp. 831-842.



\bibitem{Omel'yanov} G.A. Omel'yanov, I. Segundo-Caballero. Asymptotic and numerical description of the kink/antikink interaction. Electronic J. of Differential Equations, 2010, 150, pp. 1-19.

\bibitem{Panov} E.Yu. Panov, V.M. Shelkovich. $\delta$'-shock waves as a new type of solutions to systems of conservation laws. J. Differential Equations 228, 2006, pp. 49-86.


\bibitem{Peacock} J.A. Peacock. Cosmological Physics. 1999. Cambridge University Press.
\bibitem{Peter} P. Peter, J.Ph. Uzan. Cosmologie Primordiale. Belin, Paris, 2005. 

\bibitem{Saari} D.G. Saari, Z. Xia. Off to infinity in finite time. Notices of the AMS, 42, 5, 1995, pp. 538-546.



\bibitem{ShelkovichRMS} V.M. Shelkovich. $\delta-$ and $\delta'-$shock wave types of singular solutions of systems of conservation laws and transport and concentration processes. Russian Math. surveys 63,3, 2008, pp. 405-601.
\bibitem{Shelkovichmat} V.M. Shelkovich. The Riemann problem admitting $\delta-,\delta$'-shocks and vacuum states; the vanishing viscosity approach. J. Diff. Eq. 231, 2006, pp. 459-500.


\bibitem{Shelkovich1} V.M. Shelkovich. Transport of mass, momentum and energy in zero-pressure gas dynamics. In Proceedings of Symposia in Applied Mathematics 2009; vol.67. Hyperbolic Problems: Theory,Numerics and Applications. Edited by E. Tadmor, Jian-Guo Liu,A.E. Tzavaras. AMS, 2009, pp. 929-938.



 \bibitem{Stark} C.C. Stark, M.J. Kuchner. A new algorithm for self-consistent 3-D modeling of collisions in dusty debris disks. ArXiv.org, September 2009.
\bibitem{Temam} R. Temam. Infinite-dimensional Dynamical Systems in Mechanics and Physics. Springer Verlag, 1988. 
























  
\end{thebibliography}
 \end{document}